\newcommand{\Sy}{\mathop{\mathfrak{S}}\nolimits}
\documentclass[11pt,twoside]{amsart} 
\usepackage[OT2,T2A]{fontenc}
\usepackage{amssymb}
\usepackage{amsmath}
\advance \topmargin by -\headheight
\advance \topmargin by -\headsep
\evensidemargin \oddsidemargin
\date{}
\newtheorem{proposition}{Proposition}[section]
\newtheorem{lemma}[proposition]{Lemma}
\newtheorem{theorem}[proposition]{Theorem}
\title{Collineation group as a subgroup of the symmetric group} 
\author{Fedor Bogomolov} 
\address{F.B.:~Courant Institute of Mathematical Sciences, 251 
Mercer St., New York, NY 10012, U.S.A. \& ~Laboratory of Algebraic 
Geometry, National Research University Higher School of Economics,
7 Vavilova Str., Moscow, Russia, 117312}
\email{bogomolo@cims.nyu.edu}
\author{Marat Rovinsky} 
\address{M.R.:~Laboratory of Algebraic Geometry, 
National Research University Higher School of Economics,
7 Vavilova Str., Moscow, Russia, 117312 \& ~Institute for 
Information Transmission Problems of Russian Academy of Sciences}
\email{marat@mccme.ru}
\begin{document} 
\thanks{F.B. is supported by NSF grant DMS-1001662 
and by AG Laboratory GU-HSE grant RF government ag. 11 11.G34.31.0023}
\thanks{M.R. is supported by AG Laboratory NRU-HSE grant RF 
government ag. 11 11.G34.31.0023 and by RFBR grant 10-01-93113-CNRSL-a 
``Homological methods in geometry''}
\begin{abstract} 
Let $\Psi$ be the projectivization (i.e., the set of one-dimensional 
vector subspaces) of a vector space of dimension $\ge 3$ over a field. 
Let $H$ be a closed (in the pointwise convergence topology) subgroup 
of the permutation group $\mathfrak{S}_{\Psi}$ of the set $\Psi$. 
Suppose that $H$ contains the projective group and an arbitrary 
self-bijection of $\Psi$ transforming a triple of collinear points 
to a non-collinear triple. It is well-known from \cite{KantorMcDonough} 
that if $\Psi$ is finite then $H$ contains the alternating subgroup 
$\mathfrak{A}_{\Psi}$ of $\mathfrak{S}_{\Psi}$. 

We show in Theorem \ref{density} below that $H=\mathfrak{S}_{\Psi}$, 
if $\Psi$ is infinite. \end{abstract} 
\maketitle 

Let a group $G$ act on a set $\Psi$. For an integer $N\ge 1$, the 
$G$-action on $\Psi$ is called $N$-{\sl transitive} if $G$ acts 
transitively on the set of embeddings into $\Psi$ of a set of $N$ 
elements. This action is called {\sl highly transitive} if for any 
finite set $S$ the group $G$ acts transitively on the set of all 
embeddings of $S$ into $\Psi$. The action is highly transitive if 
and only if the image of $G$ in the permutation group 
$\mathfrak{S}_{\Psi}$ is dense in the pointwise convergence 
topology, cf. below. 

Let $\Psi$ be a projective space of dimension $\ge 2$, i.e., 
the projectivization of a vector space $V$ of dimension $\ge 3$ over 
a field $k$. The $k$-linear automorphisms of $V$ induce permutations 
of $\Psi$, called {\sl projective transformations}. 

Suppose that a group $G$ of permutations of the set $\Psi$ contains 
all projective transformations and an element which is not a 
collineation, i.e., transforming a triple of collinear points 
to a non-collinear triple. 
The main result of this paper (Theorem \ref{density}) asserts that 
under these assumptions $G$ is a dense subgroup of $\mathfrak{S}_{\Psi}$ 
if $\Psi$ is infinite. 

Somewhat similar results have already appeared in geometric context. 
We mention only some of them: \begin{itemize} \item J.Huisman and 
F.Mangolte have shown in \cite[Theorem 1.4]{HuismanMangolte} that 
the group of algebraic diffeomorphisms of a rational nonsingular 
compact connected real algebraic surface $X$ is dense in the group of 
all permutations of the set $\Psi$ of points of $X$; \item J.Koll\'ar 
and F.Mangolte have found in \cite[Theorem 1]{KollarMangolte} a 
collection of transformations generating, together with the 
orthogonal group $\mathrm{O}(3,1)$, the group of algebraic 
diffeomorphisms of the two-dimensional real sphere. \end{itemize} 

We note, however, that our result allows to work with quite 
arbitrary algebraically non-closed fields. 

{\sc Example.} Let $K|k$ be a field extension and $\tau$ be 
a self-bijection of the projective space ${\mathbb P}_k(K):=
K^{\times}/k^{\times}$, 
satisfying one of the following conditions: 
(i) $\tau:x\mapsto 1/x$ for all $x\in{\mathbb P}_k(K)$ and $K|k$ 
is separable of degree $>2$,\footnote{If $[K:k]=2$ then $\tau$ 
is projective; if $k$ is of characteristic 2 and $K\subseteq 
k(\sqrt{k})$ then $\tau$ is identical on $\mathbb{P}_k(K)$.} 
(ii) $\tau:x\mapsto x^n$ for some integer $n>1$ and all 
$x\in{\mathbb P}_k(K)$, the subfield $k$ contains all roots of 
unity of all $n$-primary degrees in $K$, and the multiplicative 
group $K^{\times}$ is $n$-divisible (e.g., if the field $K$ is 
algebraically closed). 
Then the group generated by the group $\mathrm{PGL}(K)$ of projective 
transformations ($K$ is considered as a $k$-vector space) and $\tau$ 
is $N$-transitive on ${\mathbb P}_k(K)$ for any $N$. 
Indeed, it is clear that such $\tau$'s are not 
collineations. Hence our Theorem implies the result. 

The proof of Theorem \ref{density} consists of verifying the 
$N$-transitivity of the group $H$ for all integer $N\ge 1$. 

We thank Dmitry Kaledin for encouragement and 
Ilya Karzhemanov for useful discussions and for providing 
us with a reference to the work \cite{KollarMangolte}. 
We are grateful to the referees for suggesting improvements 
of the exposition. 

\section{Permutation groups: topology and closed subgroups} 
\label{closed-subgr-symm-gr}
Let $\Psi$ be a set and $G$ be a group of its self-bijections. 
We consider $G$ as a topological group with the base of open 
subgroups formed by the pointwise stabilizers of finite 
subsets in $\Psi$. Then $G$ is a totally disconnected group. 
In particular, any open subgroup of $G$ is closed.\footnote{Indeed, 
the complement to an open subgroup is the union of translations of 
the subgroup, so it is open.} 

Denote by $\Sy_{\Psi}$ the group of all permutations of $\Psi$. 
Then the above base of open subgroups of $\Sy_{\Psi}$ is formed by 
the subgroups $\Sy_{\Psi|T}$ of permutations of $\Psi$ identical 
on $T$, where $T$ runs over all finite subsets of $\Psi$. 

\begin{lemma} \begin{enumerate} \item \label{ex-max-open} 
For any finite non-empty $T\subset\Psi$, 
$T\neq\Psi$, the normalizer $\Sy_{\Psi,T}$ of $\Sy_{\Psi|T}$ in 
$\Sy_{\Psi}$ (i.e., the group of permutations of $\Psi$, preserving 
the subset $T$) is maximal among the proper subgroups of $\Sy_{\Psi}$ 
{\rm (}\cite{F.Richman}, \cite{R.Ball}{\rm )}. 
\item Any proper open subgroup of $\Sy_{\Psi}$ which is maximal among 
the proper subgroups of $\Sy_{\Psi}$ coincides with $\Sy_{\Psi,T}$ 
for a finite non-empty $T\subset\Psi$, if $\Psi$ is infinite. 
\item Any proper open subgroup of $\Sy_{\Psi}$ is contained in 
a maximal proper subgroup of $\Sy_{\Psi}$. \end{enumerate} \end{lemma} 
{\it Proof.} By definition of our topology, any open proper subgroup $U$ 
in $\Sy_{\Psi}$ contains the subgroup $\Sy_{\Psi|T}$ for a non-empty 
finite subset $T\subset\Psi$. Assume that such $T$ is minimal. 

We claim that $\sigma(T)=T$ for all $\sigma\in U$. Indeed, if 
$\sigma(t)\notin T$ for some $t\in T$ and $\sigma\in\Sy_{\Psi}$ 
then (i) it is easy to see that the subgroup $\tilde{U}$ generated 
by $\sigma$ and $\Sy_{\Psi|T}$ meets $\Sy_{\Psi|T\smallsetminus\{t\}}$ 
by a dense subgroup, (ii) $\tilde{U}$ contains 
$\Sy_{\Psi|T\smallsetminus\{t\}}$, since both subgroups, as well as 
their intersection, are open, and thus, the intersection is closed. 
This contradicts to the minimality of $T$, and finally, $U$ is 
contained in $\Sy_{\Psi,T}$. 

The subgroups $\Sy_{\Psi,T}$ are maximal, since they are not 
embedded to each other for various $T$. \qed

\vspace{4mm}

The following `folklore' model-theoretic description of 
{\sl closed} subgroups of $\Sy_{\Psi}$ is well-known, cf. 
\cite{BeckerKechris,Cameron,Hodges}. Suppose that $\Psi$ is countable. 
Let $\mathcal{L}=\{R_i\}_{i\in I}$ be a countable relational language 
and $\mathcal{A}=(\Psi,\{R_i\}_{i\in I})$ be a structure for 
$\mathcal{L}$ with universe $\Psi$. Then $\mathrm{Aut}(\mathcal{A})$, 
the group of automorphisms of $\mathcal{A}$, is a closed subgroup of 
$\Sy_{\Psi}$. Conversely, let $H$ be a subgroup of $\Sy_{\Psi}$. 
For each $n$, let $I_n$ be the set of $H$-orbits on $\Psi^n$. 
Set $I:=\coprod_{n\ge 1}I_n$ and consider the structure 
$\mathcal{A}_H:=(\Psi,\{R_i^H\}_{i\in I})$ associated with $H$, 
where $R_i^H=i\subset\Psi^{n(i)}$. One easily checks that 
$\mathrm{Aut}(\mathcal{A}_H)$ is the closure of $H$. 

\vspace{4mm}
 
Apart from that, G.Bergman and S.Shelah prove the following 
result. Assume that $\Psi$ is countable. Let us say that two subgroups 
$G_1,G_2\subseteq\Sy_{\Psi}$ are equivalent if there exists a finite 
set $U\subseteq\Sy_{\Psi}$ such that $G_1$ and $U$ generate the same 
subgroup as $G_2$ and $U$. It is shown in \cite{BergmanShelah} that 
the closed subgroups of $\Sy_{\Psi}$ lie in precisely four equivalence 
classes under this relation. Which of these classes a closed subgroup 
$G$ belongs to depends on which of the following statements about 
open subgroups of $G$ holds: \begin{enumerate} 
\item Any open subgroup of $G$ has at least one infinite orbit in $\Psi$. 
\item There exist open subgroups $H\subset G$ such that all 
the $H$-orbits are finite, but none such that the cardinalities 
of these orbits have a common finite bound.
\item The group $G$ is not discrete, but there exist open subgroups 
$H\subset G$ such that the cardinalities of the $H$-orbits have 
a common finite bound. 
\item The group $G$ is discrete. 
\end{enumerate}
{\sc Examples of a set $\Psi$ and a group $G$ acting on it:} 
(1) (a) $\Psi$ is arbitrary and $G=\Sy_{\Psi}$. (b) $\Psi$ is a 
projective space over a field with non-discrete automorphism group 
and $G$ is the collineation group. 
(c) $\Psi$ is an infinite-dimensional projective (or affine, or 
linear) space and $G$ is the projective (or affine, or linear) group. 
(d) If $\Psi$ is arbitrary, $G$ contains a transposition $\iota$ 
of some elements $p,q\in\Psi$ and $G$ is 2-transitive then $G$ is 
a dense subgroup in $\Sy_{\Psi}$.\footnote{The transitivity does not 
suffices: if $\Psi={\mathbb Z}$ and $G$ is generated by the transposition 
$(01)$ and by the shift $n\mapsto n+2$ then the $G$-orbit of the pair 
$(0,1)$ is $\{(a,a+(-1)^a)~|~a\in{\mathbb Z}\}$, so $G$ is not dense 
in $\Sy_{\Psi}$.} (Indeed, as all transpositions generate any finite 
symmetric group, it suffices to show that $G$ contains all the 
transpositions. If $G$ is 2-transitive then for any pair of distinct 
elements $p',q'\in\Psi$ there is $g\in G$ with $g(p)=p',g(q)=q'$, and 
thus, $g\iota g^{-1}$ is the transposition the elements $p'$ and $q'$.) 
 
(4) $\Psi$ is the set of closed (or rational) points of a variety and 
$G$ is the group of points of an algebraic group acting faithfully on 
this variety; $\Psi$ is the function field of a variety and $G$ is a 
field automorphism group of $\Psi$. 

\section{Dense subgroups of the symmetric groups and the transitivity} 
It is evident that the $G$-action on $\Psi$ is highly transitive 
for any {\sl dense} subgroup $G\subseteq\Sy_{\Psi}$. 
Conversely, if a group $G$ is highly transitive on $\Psi$ then it is 
dense in $\Sy_{\Psi}$. Indeed, for any $\sigma\in\Sy_{\Psi}$ any 
neighborhood of $\sigma$ contains a subset $\sigma\Sy_{\Psi|T}$ 
for a finite subset $T\subset\Psi$, on the other hand, the identity 
embedding of $T$ into $\Psi$ and the restriction of $\sigma$ to $T$ 
belong to a common $G$-orbit, i.e., $\tau|_T=\sigma|_T$ for some 
$\tau\in G$. 

\vspace{4mm} 

We use the following terminology: (i) an $N$-{\sl set} in $\Psi$ is 
a subset in $\Psi$ of order $N$; (ii) an $N$-{\sl configuration} in 
$\Psi$ is an ordered $N$-tuple of pairwise distinct points of $\Psi$. 
The group $G$ acts naturally on the sets of $N$-sets and of 
$N$-configurations in $\Psi$. Two configurations or sets are called 
$G$-{\sl equivalent} if they belong to the same $G$-orbit. 

The $G$-action on the set of $N$-configurations in $\Psi$ commutes 
with the natural action of the symmetric group $\mathfrak{S}_N$ 
(given by $\sigma(T):=(p_{\sigma(1)},\dots,p_{\sigma(N)})$ for all 
$T=(p_1,\dots,p_N)$ and all $\sigma\in\mathfrak{S}_N$). 

\begin{lemma} For each $N\ge 1$ consider the following conditions: 
{\rm (i)}$_N$ $G$ is $N$-transitive on $\Psi$, {\rm (ii)}$_N$ 
any $N$-configuration in $\Psi$ is $G$-equivalent to a fixed 
$N$-configuration $R$ in $\Psi$, {\rm (iii)}$_N$ any 
$N$-configuration $T$ is $G$-equivalent to $\sigma(T)$ for any 
permutation $\sigma\in\mathfrak{S}_N$.\footnote{Instead of all 
permutations one can equivalently consider only a generating system 
of the group $\mathfrak{S}_N$. E.g., the transpositions (involutions 
interchanging only a pair of elements of $\{1,\dots,N\}$).} 

Then the conditions {\rm (i)}$_N$ and {\rm (ii)}$_N$ are equivalent; 
{\rm (i)}$_N$ implies {\rm (iii)}$_N$; 
{\rm (iii)}$_{N+1}$ implies {\rm (i)}$_N$ if $|\Psi|>N$. 
In particular, the group $G$ is highly transitive on 
$\Psi$ if and only if for all $N$, all $N$-configurations 
$T$ and all permutations $\sigma\in\mathfrak{S}_N$ the 
$N$-configurations $T$ and $\sigma(T)$ are $G$-equivalent. \end{lemma} 
{\it Proof.} Implications 
{\rm (iii)}$_N\Leftarrow${\rm (i)}$_N\Leftrightarrow${\rm (ii)}$_N$ 
are evident. Assume now the condition {\rm (iii)}$_{N+1}$. 

For an arbitrary pair $N$-configurations $T=(p_1,\dots,p_N)$ and 
$T'=(q_1,\dots,q_N)$ denote by $s$, $0\le s\le N$, the only integer 
such that $p_1,\dots,p_s,q_1,\dots,q_s$ are pairwise distinct and 
$p_i=q_i$ for all $i$, $s<i\le N$. To show that $T$ and $T'$ are 
$G$-equivalent, we proceed by induction on $s\ge 0$, the case $s=0$ 
being trivial. For $s>0$, the $(N+1)$-configurations $(p_1,\dots,p_N,q_1)$ 
and $(q_1,p_2,\dots,p_N,p_1)$ are $G$-equivalent, so the 
$N$-configurations $T$ and $(q_1,p_2,\dots,p_N)$ are also $G$-equivalent. 

On the other hand, the sets $\{q_1,p_2,\dots,p_N\}$ and 
$\{q_1,\dots,q_N\}$ have $N-s+1$ common elements, so 
the $N$-configurations $T'$ and $(q_1,p_2,\dots,p_N)$ 
are $G$-equivalent by the induction assumption. \qed 

\vspace{4mm} 

It is shown by H.D.Macpherson and P.M.Neumann in the usual framework 
of Zermelo--Fraenkel set theory with the axiom of choice (cf. 
\cite[Observation 6.1]{MacphersonNeumann}) that any maximal proper 
non-open subgroup of $\Sy_{\Psi}$ is dense in $\Sy_{\Psi}$. In 
particular, if $\Psi$ is a projective space then the collineation 
group (which is obviously closed) is not maximal. However, it will 
follow from Theorem \ref{density} below that the collineation group 
is maximal among proper {\sl closed} subgroups (in dimension $>1$). 

\section{Projective group as a subgroup of the symmetric group} 
In this section we prove the following 
\begin{theorem} \label{density} Let $k$ be a field, 
$V$ be a $k$-vector space of dimension $>2$ and 
$\Psi:=\mathbb{P}_k(V)=(V\smallsetminus\{0\})/k^{\times}$ 
be the projectivization of $V$. Let $H$ be a subgroup of 
$\mathfrak{S}_{\Psi}$ containing $\mathrm{PGL}(V)$ and 
an arbitrary self-bijection of $\Psi$ 
which is not a collineation. Then $H$ is dense in 
$\mathfrak{S}_{\Psi}$, if $\Psi$ is infinite. \end{theorem} 
{\it Proof.} By induction on $N$, we are going to show that $H$ is 
$N$-transitive on $\Psi$ for any $N\ge 1$. The cases $N=1,2$ are clear, 
since even the group $\mathrm{PGL}(V)$ is 2-transitive on $\Psi$. 
Though $\mathrm{PGL}(V)$ (and even the bigger group of all the 
collineations of $\Psi$) is not 3-transitive on $\Psi$ if 
$\dim_kV\ge 3$, the 3-transitivity of $H$ on $\Psi$ is evident 
(whenever $\dim_kV\ge 2$): all general 3-configurations, as well 
as all collinear 3-configurations, are $\mathrm{PGL}(V)$-equivalent, 
while any non-collineation sends some collinear 3-configuration to a 
general 3-configuration. In other words, the case $N=3$ is also trivial. 

There are two possibilities for the group $H$: 

A. There exist hyperplanes $P,P'$ in $\Psi$ and an element 
$h\in H$ such that $h(\Psi\smallsetminus P)\subseteq P'$. 
(This can happen only if $\Psi$ is infinite, i.e., either 
the field $k$ is infinite or $\Psi$ is of infinite dimension.) 

B. For any element $h\in H$ and any hyperplane $P$ in $\Psi$ the 
set $h(\Psi\smallsetminus P)$ is not contained in a hyperplane. 

\begin{lemma} \label{A+} In the case A, all $N$-configurations 
in $\Psi$ are $H$-equivalent for all $N$. \end{lemma} 
{\it Proof} proceeds by induction on $N$, the cases $N=1,2$ being 
trivial. Let $T=(p_1,\dots,p_N)$ and $T'=(p_1',\dots,p_N')$ be a 
pair of $N$-configurations. We need to show that $\xi(T)=T'$ for 
some $\xi\in H$. By the induction assumption, we may assume that 
$(p_2,\dots,p_N)=(p_2',\dots,p_N')$. 

It remains to show that for any point $p\in\Psi$ lying on neither 
of the lines passing through $p_1$ and one of the points of 
the set $\{p_2,\dots,p_N\}$ there exists an element $\xi_1\in H$ 
with $\xi_1(T)=(p,p_2,\dots,p_N)=:T''$. (Indeed, the assumption A 
implies that $|\Psi|\ge 2N|k|$, so we can choose a 
point $p$ outside the union of the $2N-2$ lines joining the 
points of the set $\{p_1',p_1\}$ and the points of the set 
$\{p_2,\dots,p_N\}$. Then the point $p_1'$ lies on neither of 
the lines passing through $p$ and one of the points of the set 
$\{p_2,\dots,p_N\}$. Therefore, $\xi_2(T'')=T'$ for some 
$\xi_2\in H$, and thus, $\xi_2\xi_1(T)=T'$.) 

By the hypothesis A, there exist hyperplanes $P,P'$ in $\Psi$ and 
an element $h\in H$ such that $h(\Psi\smallsetminus P)\subseteq P'$. 
As $h$ is surjective (and even bijective), $h(r),h(q)\notin P'$ for 
a pair of distinct points $r,q\in P$. First, we can find a projective 
transformation $g_2\in\mathrm{PGL}(V)$ sending the pair $(p_1,p)$ to 
the pair $(r,q)$ such that the support of $g_2(T)$ meets $P$ only at 
$r$. Next, we can find a projective involution $g_1\in\mathrm{PGL}(V)$ 
identical on $P'$ and interchanging the points $h(q)$ and $h(r)$. 
Then $g_2^{-1}h^{-1}g_1hg_2(T)=T''$. \qed 

\vspace{4mm}

From now on we assume the hypothesis B. 

\vspace{4mm}

Let $N\ge i\ge 1$ be integers and $T=(p_1,\dots,p_N)$ be an 
$N$-configuration $T=(p_1,\dots,p_N)$ in $\Psi$. Denote by 
$P_T^{\{i\}}$ the projective envelope of 
$p_1,\dots,\widehat{p_i},\dots,p_N$. 

We say that $T$ is $i$-{\sl disjoint} if $p_i\notin P_T^{\{i\}}$, 
and that $T$ is {\sl disjoint} if $T$ is $i$-disjoint for some $i$. 

\begin{lemma} \label{extra-pt-triv-B-i} Let $N\ge i\ge 1$ be 
integers. Assume that $H$ is $(N-1)$-transitive on $\Psi$. 

Then, in the case B, all $i$-disjoint 
$N$-configurations are $H$-equivalent. \end{lemma} 
{\it Proof.} Let $T=(p_1,\dots,p_N)$ and $T'=(p_1',\dots,p_N')$ 
be $i$-disjoint $N$-configurations for some $i$ 
(so $\dim\Psi\ge 2$ if $N\ge 3$). 

As $H$ is $(N-1)$-transitive on $\Psi$, we can choose an element 
$h\in H$ with $h(p_1,\dots,\widehat{p_i},\dots,p_N)=
(p_1',\dots,\widehat{p_i'},\dots,p_N')$. As $h(\Psi\smallsetminus 
P_T^{\{i\}})$ is not contained in $P_{T'}^{\{i\}}$, there exist: 
(i) a point $p\in\Psi\smallsetminus P_T^{\{i\}}$ such 
that $hp\not\in P_{T'}^{\{i\}}$, (ii) $g_2\in\mathrm{PGL}(V)$ 
identical on $P_T^{\{i\}}$ and sending $p_i$ to $p$, (iii) 
$g_1\in\mathrm{PGL}(V)$ identical on $P_{T'}^{\{i\}}$ and sending 
$h(p)$ to $p_i'$. Therefore, $g_1hg_2(T)=T'$, as desired. \qed

\begin{lemma} \label{extra-pt-triv-B} Assume that $H$ 
is $(N-1)$-transitive on $\Psi$ for an integer $N\ge 1$. 

Then, in the case B, all disjoint $N$-configurations are 
$H$-equivalent. In particular, the permutation group 
$\mathfrak{S}_N$ preserves the $H$-equivalence class 
of any disjoint $N$-configuration. \end{lemma} 
{\it Proof.} Fix some pair $i\neq j$ with $1\le i,j\le N$. Let us 
show that any $i$-disjoint $N$-configuration $T=(p_1,\dots,p_N)$ in 
$\Psi$ is $H$-equivalent to a $j$-disjoint $N$-configuration. As $H$ 
is $(N-1)$-transitive on $\Psi$, we can choose an element $\xi\in H$ 
such that $\xi(p_s)=p_s$ for all $s\neq i,j$ and $\xi(p_j)=p_i$. 
If $\xi(p_i)\in P_T^{\{i\}}$ then $\xi(T)$ is $j$-disjoint. If 
$\xi(p_i)\notin P_T^{\{i\}}$ then there exists a projective 
involution $\iota$ fixing $P_T^{\{i\}}$ and interchanging 
$\xi(p_i)$ and $p_i$, i.e., $\xi^{-1}\iota\xi(T)$ (having 
the same support as $T$) is $j$-disjoint. In both cases $T$ 
is $H$-equivalent to a $j$-disjoint $N$-configuration, while 
all $j$-disjoint $N$-configurations are $H$-equivalent. \qed 

\vspace{4mm} 

The following Lemma reduces verification of the $N$-transitivity 
to checking of the $H$-equivalence of all $N$-sets. 
\begin{lemma} 
Let $P$ be a hyperplane in $\Psi$. Suppose that $H$ is 
$(N-1)$-transitive on $\Psi$. 

Then, in the case B, the permutation group $\mathfrak{S}_N$ 
preserves the $H$-equivalence class of any $N$-configuration in 
$\Psi$, if $|P|\ge N-2$. \end{lemma} 
{\it Proof.} Suppose that an $N$-configuration $T$ in $\Psi$ 
is not $H$-equivalent to a disjoint one. By $(N-1)$-transitivity 
of $H$, we may assume that $T=(q_1,\dots,q_N)$, where 
$(q_1,\dots,\widehat{q_i},\dots,\widehat{q_j},\dots,q_N)$ is 
a fixed $(N-2)$-configuration in $P$ for some pair $i\neq j$ 
with $1\le i,j\le N$, and $q_i\in\Psi$ is a fixed point outside 
of $P$. Then, as $T$ is not $H$-equivalent to a disjoint 
configuration, $q_j\in\Psi\smallsetminus P$. 

To show that $T$ is $H$-equivalent to 
$\sigma(T):=(q_{\sigma(1)},\dots,q_{\sigma(N)})$ for any permutation 
$\sigma\in\mathfrak{S}_N$ it suffices to verify that $T$ is 
$H$-equivalent to $\sigma_{ij}(T)$ for the transposition $\sigma_{ij}$ 
of any pair $1\le i<j\le N$. But this is clear, since there exists 
a projective involution $\iota_{ij}$ fixing $P$ and 
interchanging $q_i$ and $q_j$. \qed 

\vspace{4mm}

For any pair of subsets $\Pi_1,\Pi_2\subset\Psi$ we introduce the subset 
$H_{\Pi_1,\Pi_2}:=\{h\in H~|~h(\Pi_1)\subseteq\Pi_2\}$ in $H$. There is 
a natural composition law: $H_{\Pi_1,\Pi_2}\times H_{\Pi_2,\Pi_3}
\to H_{\Pi_1,\Pi_3}$. 
One has $gH_{\Pi_1,\Pi_2}=H_{\Pi_1,g\Pi_2}$ and 
$H_{\Pi_1,\Pi_2}g=H_{g^{-1}\Pi_1,\Pi_2}$ for any $g\in H$. 

Define a binary relation $\succ_{\Pi_1,\Pi_2}$ on $\Psi$ by 
the condition $x\succ_{\Pi_1,\Pi_2}y$ if and only if one 
has $h(y)\in\Pi_2$ for any $h\in H_{\Pi_1,\Pi_2}$ such that 
$h(x)\in\Pi_2$. Clearly, $g(x)\succ_{g(\Pi_1),\Pi_2}g(y)$ 
for any $g\in H$ if $x\succ_{\Pi_1,\Pi_2}y$. 
\begin{lemma} \label{equiv-rel} The binary relation 
$\succ_{\Pi_1,\Pi_2}$ is reflexive and transitive. Moreover, if 
$\Pi_1$ is contained in a proper projective subspace $P$ in $\Psi$ 
then {\rm (i)} the restriction of $\succ_{\Pi_1,\Pi_2}$ to 
$\Psi\smallsetminus P$ is an equivalence relation; {\rm (ii)} 
the equivalence classes in $\Psi\smallsetminus P$ are complements 
to $P$ of projective subspaces in $\Psi$. \end{lemma} 
{\it Proof.} The reflexivity and the transitivity of 
$\succ_{\Pi_1,\Pi_2}$ are trivial. 
(i) If $\Pi_1$ is contained in a subspace $P\subset\Psi$ then for 
any pair of points $x,y\in\Psi\smallsetminus P$ there is a projective 
involution $\iota$ identical on $P$ and interchanging $x$ and $y$. 
Now, if $x\succ_{\Pi_1,\Pi_2}y$ then $h(y)\in\Pi_2$ for any 
$h\in H_{\Pi_1,\Pi_2}$ such that $h(x)\in\Pi_2$. As 
$H_{\Pi_1,\Pi_2}=H_{\Pi_1,\Pi_2}\iota$, one has $h'\iota(y)\in\Pi_2$ 
for any $h'\in H_{\Pi_1,\Pi_2}$ such that $h'\iota(x)\in\Pi_2$, 
i.e., $y\succ_{\Pi_1,\Pi_2}x$. 

(ii) A subset in $\Psi\smallsetminus P$ is a complement to $P$ 
of a projective subspace in $\Psi$ if and only if together 
with a pair of points $x,y\in\Psi\smallsetminus P$ it contains 
the line $\overline{xy}$ passing through them (eventually, 
punctured at the meeting point with $P$). Thus, we need to 
show that for any triple of pairwise distinct collinear points 
$x,y,z\in\Psi\smallsetminus P$ such that $x\succ_{\Pi_1,\Pi_2}y$ 
one has $x\succ_{\Pi_1,\Pi_2}z$. Indeed, there is a projective 
transformation $\alpha$ identical on $P\cup\{x\}$ and sending 
$y$ to $z$. As $H_{\Pi_1,\Pi_2}=H_{\Pi_1,\Pi_2}\alpha$, one 
has $h'\alpha(y)\in\Pi_2$ for any $h'\in H_{\Pi_1,\Pi_2}$ such 
that $h'\alpha(x)\in\Pi_2$, i.e., $x\succ_{\Pi_1,\Pi_2}z$. \qed

\begin{lemma} \label{gen-nongen-extra-pt-triv-B} Suppose that 
$H$ is $(N-1)$-transitive on $\Psi$. Then, in the case B, any 
$N$-set is $H$-equivalent to a disjoint one,\footnote{Naturally, 
a set is called {\sl disjoint} if one of its points is not in the 
projective envelope of the others.} if $N\le|k|+2$. \end{lemma} 
{\it Proof.} By $(N-1)$-transitivity of $H$, if 
$N\le|k|+2=\#{\mathbb P}^1(k)+1$ then any $N$-set in $\Psi$ 
is $H$-equivalent to $T=\{q_1,q_2\}\cup R$, where $R$ is a 
fixed $(N-2)$-subset of a projective line $l\subset\Psi$ and 
$q_2\in\Psi$ is a fixed point. If $N=|k|+2=\#{\mathbb P}^1(k)+1$ 
then taking $q_2\in l$ we get a disjoint $T$. 

From now on $N\le|k|+1$. Then fixing $q_2\in\Psi\smallsetminus l$ 
we may assume that $q_1\in\Psi\smallsetminus l$, that 
$q_2\succ_{R,l}q_1$ and that $T$ is coplanar, as otherwise $T$ 
is $H$-equivalent to a disjoint $N$-set. 

Let $q_0$ be the intersection point of $l$ and the line 
$\overline{q_1q_2}$ passing through $q_1$ and $q_2$. 

Suppose first that $q_0\notin R$. Then $\{q_0,q_2\}\cup R$ 
is a disjoint $N$-set and Lemma \ref{extra-pt-triv-B} implies 
that there is $\xi\in H$ interchanging $q_0$ and $q_2$ and 
identical on $R$. Then $\xi(q_0)\notin l$. By Lemma 
\ref{equiv-rel} (ii), where we take $\Pi_1=R$ and $\Pi_2=l$, 
one has $\xi(\overline{q_1q_2}\smallsetminus\{q_0\})\subset l$. 
We can choose such $\psi\in H$ that 
$\psi(R)\subset\overline{q_1q_2}$ and $\psi(q_2)=q_0$. 
If $\psi(T)$ is not a disjoint set then 
$\psi(q_1)\in\overline{q_1q_2}$, but 
then $\xi\psi(T)$ is a disjoint set. 

This settles the case $q_0\notin R$, so let us now suppose 
that $q_0\in R$ and that $T$ is not $H$-equivalent 
to a disjoint $N$-set. If we change $q_2$ in its 
$\succ_{R,l}$-equivalence class then the resulting new $T$ is 
not $H$-equivalent to a disjoint $N$-set. Then 
$\overline{q_1q_2}\smallsetminus\{q_0\}$ is precisely the 
$\succ_{R,l}$-equivalence class of $q_1$. (Otherwise, by 
Lemma \ref{equiv-rel} (ii), we can choose a new $q_2$ in its 
equivalence class so that the new $q_0$ is not in $R$, 
but then $T$ is $H$-equivalent to a disjoint $N$-set.) 
We can find an element $\gamma\in H$ such that all 
the points in the support of $\gamma(T)$ are on the line 
$\overline{q_1q_2}$, and $\gamma(q_1)=q_0\in R\subset l$. 

Then any $\beta\in H$ fixing $\gamma(q_2)$ and inducing 
a cyclic permutation of $R$ induces an automorphism of 
the $\succ_{R,l}$-equivalence class 
$\overline{q_1q_2}\smallsetminus\{q_0\}$, 
so $\beta\gamma(T)$ is disjoint. \qed 

\section{Case of finite field $k$ and the end of the proof}
\begin{lemma} \label{choice-of-B} 
Let $V$ be a vector space over a finite field $k$, $\tilde{h}$ be a 
self-embedding of $V\smallsetminus\{0\}$, $J\subset V$ be a finite 
set such that $J$ and $\tilde{h}(J)$ consist of independent vectors 
and $J\subset P_0\subset P_1\subset P_2\subset\dots\subset V=
\bigcup_iP_i$, $\dim P_i=|J|+i$, be a flag of vector subspaces. 
Then there exists a basis $\mathcal{B}=J\cup\{e_1,e_2,\dots\}$ 
of $V$ such that $e_i\in P_i\smallsetminus P_{i-1}$ and 
$\tilde{h}(\mathcal{B})$ consists of independent vectors. \end{lemma} 
{\it Proof.} It is possible to choose such $e_i$ inductively, since 
$\#\tilde{h}(P_i\smallsetminus P_{i-1})=\#(P_i\smallsetminus P_{i-1})
=(\# k-1)(\# k)^{|J|+i-1}>\#[\langle\tilde{h}(J),\tilde{h}(e_1),
\dots,\tilde{h}(e_{i-1})\rangle\smallsetminus\{0\}]=\#(P'_{i-1}
\smallsetminus\{0\})=(\# k)^{|J|+i-1}-1$. \qed

\begin{lemma} \label{gen-nongen-B-finite-field} 
Suppose that $H$ is $(N-1)$-transitive on $\Psi$, the field 
$k\cong{\mathbb F}_q$ is finite of order $q$ and either $\Psi$ 
is infinite-dimensional, or $\dim\Psi\ge N-1\ge 3$ and $q>2$. 
Then, in the case B, any $N$-set in $\Psi$ is $H$-equivalent 
to a general $N$-set in $\Psi$. \end{lemma} 
{\it Proof.} For any subset $A\subset\Psi$ denote by $P_A$ the 
projective envelope of $A$. Let $P$ be an $(N-2)$-dimensional 
subspace in $\Psi$. There is a disjoint $N$-set in $P$: 
$N\le\#{\mathbb P}^{N-3}({\mathbb F}_q)+1=\frac{q^{N-2}-1}{q-1}+1$. 
Then, by Lemma \ref{extra-pt-triv-B}, there exists $h\in H$ such 
that $h(P)$ is contained in no $(N-2)$-dimensional subspace. 
Fix such $h$. Let $S$ be a maximal independent subset in $P$ 
such that $h(S)$ is also independent. 

By Lemma \ref{choice-of-B}, one has $P_S=P$ (in other words, $|S|=N-1$). 

By $(N-1)$-transitivity, any $N$-set is $H$-equivalent to an $N$-set 
$T=\{p_1\}\cup S$ in $\Psi$. If either $p_1\in P$ and $p_1$ is not in 
general position with respect to $S$, or $p_1\notin P$ then $T$ is 
disjoint, so by Lemma \ref{extra-pt-triv-B}, $T$ is $H$-equivalent to a 
general $N$-set. Suppose therefore that $T$ is a general $N$-set in $P$. 

As all general $N$-sets in $P$ are $\mathrm{PGL}(V)$-equivalent, 
we may assume that $h(p')\in P_{h(S)}$ and that $h(p')$ is in 
general position with respect to $h(S)$ for any point $p'\in P$ 
in general position with respect to $S$. Note, however, that 
there is only one such $p'$ in the case $q=2$. 

Fix some $p\in P$ such that $h(p)\notin P_{h(S)}$. In particular, 
the set $\{p\}\cup S$ is $H$-equivalent to a general set. 
Then $p$ is a point of $P_I$ in general position with respect 
to $I$ for some subset $I\subsetneqq S$ with $|I|\ge 2$. 

Fix some $s\in I$ and choose homogeneous coordinates 
$X_2,\dots,X_N$ on $P$ such that the elements of $S$ are given 
by $X_2=\dots=\widehat{X_i}=\dots=X_N=0$ for $2\le i\le N$, 
the elements of $I$ correspond to $2\le i\le|I|+1$, the 
point $s$ corresponds to $i=|I|+1$ and the point $p$ is given 
by $X_2=\dots=X_{|I|+1}\neq 0$ $\&$ 
$X_{|I|+2}=\dots=\widehat{X_i}=\dots=X_N=0$. 

In the case $q=2$, our choice of $h$ contradicts to the 
conclusion of Lemma \ref{quasi-proj} below, and thus, to the 
assumption that $T$ is not $H$-equivalent to a general set.

In the case of $q>2$, we are looking for a point $s'\in P$ 
in general position with respect to $S$ and in general 
position with respect to $\{p\}\cup(S\smallsetminus\{s\})$. 
In coordinates: $X_2\cdots X_N\neq 0$ and 
$\prod_{i=1}^{|I|}(X_i-X_{|I|+1})\neq 0$. 
As $q\ge 3$, we can find such a point. 

Then $\{p,s'\}\cup S\smallsetminus\{s\}$ is a general $N$-set 
in $P$ and its image under $h$ is disjoint, i.e., it is 
$H$-equivalent to a general $N$-set. \qed

\begin{lemma} \label{quasi-proj} Let $m,n\ge 1$ be integers, and $P$ 
be an $m$-dimensional vector subspace of an infinite-dimensional 
${\mathbb F}_2$-vector space $V$. Let $h$ be a self-bijection of $V$ 
such that $h(0)=0$. Suppose that the image under $h$ of any general 
$(n+1)$-set in any $n$-dimensional vector subspace is not general. 

Then $h(P)$ is a vector subspace of $V$ isomorphic to $P$. \end{lemma} 
{\it Proof.} Set $m=\dim P$ and consider an $(m+1)$-dimensional 
vector subspace $\tilde{P}$ in $V$ containing $P$. 
By Lemma \ref{choice-of-B}, there exist a basis $\mathcal{B}$ of 
$V$ such that (i) the vectors of $\tilde{h}(\mathcal{B})$ are 
independent, (ii) $\mathcal{B}$ contains bases of $P$ and of $\tilde{P}$. 
If $x$ is the unique element of $\mathcal{B}\cap(\tilde{P}\smallsetminus 
P)$, we write also $\mathcal{B}=\mathcal{B}_x$. 

Denote by $V^+_x$ the hyperplane in $V$ of sums of even number of 
elements of $\mathcal{B}_x$. Denote by $P^+$ (resp., by $\tilde{P}^+_x$) 
the corresponding hyperplane in $P$ (resp., in $\tilde{P}$, 
$\tilde{P}^+_x$ does not contain $x$). There are precisely 
three hyperplanes in $\tilde{P}$ containing $P^+$: (i) $P$, 
(ii) a hyperplane containing $x$, (iii) $\tilde{P}^+_x$. 

Set $V^-_x:=P\smallsetminus V^+_x$, $P^-:=P\smallsetminus P^+$ and 
$\tilde{P}^-_x:=\tilde{P}\smallsetminus\tilde{P}^+_x$. Suppose that $h$ 
does not transform a general $(n+1)$-set in any $n$-dimensional vector 
subspace to a general one. Then for any general set $I\subset V$ of order 
$n$ with general $h(I)$ one has $h(\sum_{v\in I}v)=\sum_{v\in I}h(v)$. 
By Lemma \ref{bool-alge}, $h(V^-_x)$ is an affine subspace of the span 
of $h(\mathcal{B}_x)$ and $h|_{V^-_x}$ is the restriction of a linear 
endomorphism of $V$. 

As $P\cap\tilde{P}^+_x=P^+$, one has 
$\#(\tilde{P}\smallsetminus(P\cup\tilde{P}^+_x))=2^m-2^{m-1}=2^{m-1}
>\#(P\smallsetminus P^-)=2^{m-1}-1$, and thus, there is 
$y\in\tilde{P}\smallsetminus(P\cup\tilde{P}^+_x)$ with $h(y)\notin P$. 
Therefore, $\tilde{P}=P\cup\tilde{P}^+_x\cup\tilde{P}^+_y$. 

Then the restriction of $h$ to $\tilde{P}^-_x=P^-\cup(\tilde{P}^+_y
\smallsetminus P^+)$ coincides with the restriction of a linear map, 
and hence, if $t\in\tilde{P}^-_x$ and $h(t)$ is in the linear span $L$ 
of $h(P\cap\mathcal{B}_x)$ then $t$ is in the linear envelope of $P^-$, 
i.e., $t\in P$. Similarly, if $t\in\tilde{P}^-_y$ and $h(t)$ is in $L$ 
then $t\in P$. As $\tilde{P}=P\cup\tilde{P}^-_x\cup\tilde{P}^-_y$, we 
get that $h(x_0)$ is in $L$ a point $x_0\in\Psi$ only if $x_0\in P$. 
We conclude that $h(P)$ contains $L$, since $h$ is surjective (and even 
bijective), and thus, $h(P)=L$, since $P$ is finite and $L\cong P$. \qed 

\begin{lemma} \label{bool-alge} 
Let $\mathcal{B}$ be a basis in an infinite-dimensional 
${\mathbb F}_2$-vector $V$ and $n\ge 2$ be an integer. Denote by 
$V^+$ the hyperplane in $V$ consisting of all sums of even number 
of elements of $\mathcal{B}$. A subset of $V$ is called 
$n$-{\rm closed} if it contains the sums of all 
collections of its $n$ independent elements. 
\begin{enumerate} \item \label{n-span} Let $\mathcal{C}$ be 
the minimal $n$-closed subset of $V$ containing all one-element 
subsets of $\mathcal{B}$. Then $\mathcal{C}=V\smallsetminus\{0\}$ 
if $n$ is even, $\mathcal{C}$ consists of all sums of odd number 
of elements of $\mathcal{B}$ if $n$ is odd.\footnote{This is an 
affine space over $V^+$.} \item \label{one-all} Let 
$\tilde{\mathcal{C}}$ be the minimal $n$-closed subset of 
$V$ containing all sums of odd number of elements of 
$\mathcal{B}$ and a non-zero vector $I\in V^+$. Then 
$\tilde{\mathcal{C}}=V\smallsetminus\{0\}$. \end{enumerate} \end{lemma} 
{\it Proof.} Note, the map $J\mapsto\sum_{v\in J}v$ induces an 
isomorphism of the group of the finite subsets in $\mathcal{B}$ 
(with the operation $\Delta$ of symmetric difference) onto $V$. 
We have to show that $\mathcal{C}$ contains the sum of all elements 
of an $m$-subset $J\subset\mathcal{B}$ for any $m\ge 1$, odd in 
the case of odd $n$. Such a $J$ is constructed as 
$I_1\Delta\dots\Delta I_n$ for some independent 
sets $I_1,\dots,I_n\in\mathcal{C}$, uniquely (modulo 
$\mathfrak{S}_{\mathcal{B}}$-action) determined by the numbers 
$|I_1|$, $|I_1\cap(I_2\cup\dots\cup I_n)|=a$ and some extra conditions 
described below. For any $I_1,\dots,I_n\in\mathcal{C}$ one has 
$|J|\equiv\sum_{j=1}^n|I_j|\pmod 2$. In particular, 
$|J|\equiv n\equiv 1\pmod 2$ if $n$ is odd. In all cases we impose 
the conditions $|I_2|=\dots=|I_n|=1$ and $|I_1|+1\equiv n+m\pmod 2$. 

If $m\le 2n-1$ is odd, it suffices to impose the conditions 
$|I_1|=n$, $a=n-\frac{m+1}{2}$. 

If $m\le 2n-2$ is even and $n$ is even, it suffices to 
ask that $|I_1|=n-1$ (is odd, so $I_1\in\mathcal{C}$), 
$a=n-1-\frac{m}{2}$.

If $m\ge n$ ($m$ is odd, if $n$ is odd), we proceed by induction 
on (odd, if $n$ is odd) $m$ and take some disjoint sets 
$I_1,\dots,I_n\in\mathcal{C}$, where $|I_1|=m-n+1$. 

Suppose that $n$ is odd. To show that $\tilde{\mathcal{C}}$ contains 
the sum of all elements of an $m$-set for any even $m\ge 2$, we take 
some $(m+|I|+n-2)$-set $J$ containing $I$ and choose pairwise 
distinct one-element subsets $J_3,\dots,J_n\subset J\smallsetminus I$. 
Then $|J\Delta I\Delta J_3\Delta \dots\Delta J_n|=m$. \qed 

\vspace{4mm} 

The proof of Theorem \ref{density} is concluded by combining 
Lemma \ref{extra-pt-triv-B} either with Lemma 
\ref{gen-nongen-extra-pt-triv-B} (in the case of infinite $k$), or with 
Lemma \ref{gen-nongen-B-finite-field} (in the case of finite $k$). \qed

\section{Remarks on the case of finite $\Psi$} 
\begin{theorem}[Kantor--McDonough, \cite{KantorMcDonough}] 
\label{density-ii} Let $k$ be a finite field 
and $V$ be a $k$-vector space of a finite dimension $>2$. Set 
$\Psi:=\mathbb{P}_k(V)=(V\smallsetminus\{0\})/k^{\times}$. Let $H$ 
be a subgroup of $\mathfrak{S}_{\Psi}$ containing $\mathrm{PGL}(V)$ 
and an arbitrary self-bijection of $\Psi$ which is not a collineation. 
Then $H$ contains the alternating subgroup $\mathfrak{A}_{\Psi}$. 
\end{theorem} 

For the sake of completeness we mention some details of the 
proof of Theorem \ref{density-ii}. 
Let $G$ be a group acting on a projective space 
${\mathbb P}^d({\mathbb F}_q)$. Suppose that $G$ contains the 
special projective group, but its action is not $m$-transitive for 
$m=\#{\mathbb P}^{d-1}({\mathbb F}_q)$ then 
\cite[Theorem 1.1 (iii)]{Kantor} gives a choice of 3 groups of 
possibilities for $G$, one is being excluded by the principal theorem 
of projective geometry, while the remaining 2 being excluded for 
arithmetical reasons. Then it remains to prove the following 
\begin{proposition} 
Let $d\ge 2$ be an integer and ${\mathbb F}_q$ be a finite field. 
Suppose that a group $G$ acts $m$-transitively 
on the projective space ${\mathbb P}^d({\mathbb F}_q)$, 
where $m=\#{\mathbb P}^{d-1}({\mathbb F}_q)$. 
Then $G$ contains the alternating group 
$\mathfrak{A}_{{\mathbb P}^d({\mathbb F}_q)}$. 
\end{proposition} 
{\it Proof.} A theorem of H.Wielandt from \cite{Wielandt} asserts 
that {\it if a group $G$ permuting $M$ elements is $m$-transitive 
then $m<3\log(M-m)$, unless $G$ contains $\mathfrak{A}_M$}. 

The following calculation lists all 9 cases where the condition 
$m<3\log(M-m)$ of Wielandt's theorem is not satisfied. 
\begin{lemma} \label{appl-Wielandt} 
Let $d\ge 2$ be an integer and ${\mathbb F}_q$ be a finite field. 
Suppose that $\#{\mathbb P}^{d-1}({\mathbb F}_q)<
3\log\#{\mathbb A}^d({\mathbb F}_q)$. Then $d=2$ and $q\le 13$. \end{lemma} 
{\it Proof.} Consider the function $\xi(d)=\#{\mathbb P}^{d-1}
({\mathbb F}_q)-3\log\#{\mathbb A}^d({\mathbb F}_q)
=\frac{q^d-1}{q-1}-3d\log q$. Its derivative 
$\xi'(d)=\frac{q^d\log q}{q-1}-3\log q$ vanishes only at 
$\frac{\log(3(q-1))}{\log q}$. We see from $(q-3/2)^2+3/4>0$ 
that $3(q-1)<q^2$, and thus, $\frac{\log(3(q-1))}{\log q}<2$. In other 
words, $\xi$ increases in $d\ge 2$. 

Now consider the (same) function $\varphi(q)=
\#{\mathbb P}^{d-1}({\mathbb F}_q)-3\log\#{\mathbb A}^d({\mathbb F}_q)
=\frac{q^d-1}{q-1}-3d\log q$. 

If $d=3$ then $\varphi(q)=q^2+q+1-9\log q$. As $\varphi'(q)=2q+1-9/q$, 
the critical point of $\varphi$ is at $\frac{\sqrt{73}-1}{4}<2$, so 
$\varphi$ increases in $q\ge 2$. As $\varphi(2)=7-9\log 2>0$, the 
function $\varphi$ is positive. 

If $d=2$ then $\varphi(q)=q+1-6\log q$. As $\varphi'(q)=1-6/q$, 
the critical point of $\varphi$ is at $6$. One has 
$\varphi(2)=3-6\log 2<0$, $\varphi(15)=16-6\log 15<0$, 
$\varphi(16)=17-6\log 16>0.3>0$, so the function 
$\varphi$ is positive for $q\ge 16$. This means that $\varphi(q)$ 
is negative only if $q$ is one of $2,3,4,5,7,8,9,11,13$. \qed 

\vspace{4mm}

In the cases $q\neq 2,4$ we apply the following theorem 
of G.Miller from \cite{Miller}: 
{\it If $M=qp+m$, where $p$ is a prime, $p>q>1$ and $m>q$, then 
a group permuting $M$ elements can be at most $m$-transitive, 
unless it contains the alternating group.} 
Namely, setting $m_0$ for a lower bound of transitivity, 
$$\begin{array}{cccc}
q&M&m_0&M=qp+m\\
3&13&4&13=2\cdot 5+3\\
5&31&6&31=2\cdot 13+5\\
7&57&8&57=3\cdot 17+6\\
8&73&9&73=5\cdot 13+8\\
9&91&10&91=2\cdot 41+9\\
11&133&12&133=2\cdot 61+11\\
13&183&14&183=10\cdot 17+13
\end{array}$$ 
\vspace{4mm}

In the case $q=2$, i.e., of the projective plane over 
$k={\mathbb F}_2$ (with $(M,m_0)=(7,3)$), one has 
$\#\mathrm{PGL}_3({\mathbb F}_2)=
7\cdot 6\cdot 4$, so $[\mathfrak{A}_7:\mathrm{PGL}_3({\mathbb F}_2)]
=5\cdot 3=15$. By Lemma \ref{gen-nongen-extra-pt-triv-B}, $G$ is 
$3$-transitive. There are precisely $\binom{7}{3}=35$ $3$-sets, 
so $\#G$ is divisible by the least common multiple 
$7\cdot 6\cdot 5\cdot 4$ of $7\cdot 6\cdot 4$ and $35$, 
and thus, $[\mathfrak{S}_7:G]|6$. 
The only such subgroups are $\mathfrak{A}_7$ and $\mathfrak{S}_7$.

\vspace{4mm}

In the remaining case $q=4$ (with $(M,m_0)=(21,5)$) we apply 
the following theorem of C.Jordan: 
{\it If a primitive group contains a cycle of length $p$ 
and permutes $M=p+m$ elements, where $p$ is a 
prime and $m>2$, then it contains the alternating group.} 

\vspace{4mm}

{\it Remarks.} 1. The collineation group of $\mathbb{P}_k(V)$ is 
maximal among proper {\sl closed} subgroups of $\mathfrak{S}_{\Psi}$. 
 
2. {\sl Parity of projective transformations.} 
Let $\mathbb{F}_q$ be a finite field and $n\ge 1$ be 
an integer. The projective special linear group 
$\mathrm{PSL}_{n+1}(\mathbb{F}_q)$ is simple, with two solvable 
exceptions: $\mathrm{PSL}_2(\mathbb{F}_2)\cong\mathfrak{S}_3$ 
and $\mathrm{PSL}_2(\mathbb{F}_3)\cong\mathfrak{A}_4$. 

In any non-exceptional case, any element of the maximal abelian 
quotient of the projective group $\mathrm{PGL}_{n+1}(\mathbb{F}_q)$ 
is presented by a diagonal matrix $g_{\lambda}:=\mathrm{diag}
(\lambda,1\dots,1)\in\mathrm{PGL}_{n+1}(\mathbb{F}_q)$ for some 
$\lambda\in\mathbb{F}_q^{\times}$. Let $s\ge 1$ be minimal with 
$\lambda^s=1$. Then $g_{\lambda}$ acts on an $n$-dimensional 
projective space over $\mathbb{F}_q$ with a fixed hyperplane 
and an extra fixed point. Other orbits consist of $s$ elements, 
and therefore, the parity of $g_{\lambda}$ coincides with the 
parity of $\frac{q^n-1}{s}(s-1)=q^n-1-\frac{q^n-1}{s}\equiv 
q+1-\frac{q-1}{s}(1+(n-1)q)\pmod 2$. 
Finally, $\mathrm{PGL}_{n+1}(\mathbb{F}_q)$ contains 
an odd permutation if and only if $qn$ is odd. 

3. Let $G$ be a finite permutation group, which is neither symmetric 
nor alternating group. As mentioned in \cite[\S7.3]{DixonMortimer}, 
it can be deduced from the classification of finite simple groups that 
$G$ is at most 5-transitive; moreover, if $G$ is 4- or 5-transitive 
then $G$ is one of the Mathieu groups $M_{11}$, $M_{12}$, $M_{23}$, 
$M_{24}$. 

\end{document}